\documentclass[12pt]{article}
\usepackage{amsfonts}
\usepackage{stmaryrd}
\usepackage{pifont}
\usepackage{amssymb,amsmath,amsfonts,mathrsfs,amscd}
\usepackage{shadow}
\usepackage{flushend, cuted}
\usepackage{multicol}
\usepackage{color}
\usepackage{longtable,colordvi}
\newtheorem{thm}{Theorem}[section]
\newtheorem{prop}[thm]{Proposition}

\newtheorem{lem}[thm]{Lemma}

\def\qed{\nopagebreak\hfill{\rule{4pt}{7pt}}\medbreak}

\makeatletter \@addtoreset{equation}{section} \makeatother
\begin{document}
\begin{center}
{\Large\bf Eulerian Pairs on Fibonacci Words}
\end{center}

\begin{center}
Teresa X.S. Li$^1$, Charles B. Mei$^2$, Melissa Y.F. Miao$^3$\\[6pt]
$^{1}$School of Mathematics and Statistics
Southwest University\\[6pt]
Chongqing 400715, P.R. China\\[6pt]
$^{2,3}$Center for Combinatorics, LPMC-TJKLC\\[6pt]
Nankai University, Tianjin 300071, P.R. China\\[6pt]
$^1$pmgb@swu.edu.cn, $^2$meib@mail.nankai.edu.cn,
$^3$miaoyinfeng@mail.nankai.edu.cn
\end{center}

\begin{abstract}
 Recently, Sagan and Savage introduced the
 notion of Eulerian pairs. In this note, we find Eulerian pairs on Fibonacci words based on Foata's first transformation or Han's bijection and a map
 in the spirit of a bijection of Steingr\'{\i}msson.

\end{abstract}



\section{Introduction}

This paper is motivated by the notion of Eulerian pairs introduced
by Sagan and Savage \cite{Sagan-Savage2011} in their study of
Mahonian pairs. Let $\mathbb{P}$ be the set of positive integers and
let $\mathbb{P}^{*}$ be the set of words on $\mathbb{P}$. For two
finite subsets $S,T \subset \mathbb{P}^{*}$, the pair $(S,T)$ is
called a Mahonian pair if the distribution of the major index over
$S$ is the same as the distribution of the inversion number over
$T$.  Similarly, $(S,T)$ is said to be an  Eulerian pair if the
distribution of the descent number over $S$ is the same as the
distribution of the excedance number over $T$.

 The well-known theorem of MacMahon  \cite{MacMahon1916} can be rephrased as the fact that $(S_n,S_n)$ is a Mahonian pair, where $S_n$ is the set of
permutations on $[n]=\{1,2,\ldots,n\}$.
 Foata \cite{Foata1968} found a combinatorial proof of this fact
  by establishing a correspondence which has been called
    the second fundamental transformation,
denoted $\Phi_2$. With the aid  of the map $\Phi_2$, Sagan and Savage found
 Mahonian pairs $(S,\Phi_2(S))$, where $S$ is a set of ballot sequences or a set of
 Fibonacci words. By a Fibonacci word we mean a word on
$\{1,2\}$ containing no consecutive ones. Dokes et al.\cite{Dokos} studied mahonian pairs on permutations avoiding some patterns.
 In this paper, we find Eulerian pairs on Fibonacci words based on bijections of
  Foata \cite{Foata1965}, Han\cite{Han1994} and Steingr\'{\i}msson \cite{Steingrimsson1991}.

  We  adopt some common notation on words.
 For a  word $\omega=a_1a_2\cdots a_n$ ,
  the descent number $\mathrm{des}(\omega)$,  the inversion number $\mathrm{inv}(\omega)$ and the major index $\mathrm{maj}(\omega)$ are defined by
\begin{align*}
&\mathrm{des}(\omega)=\#\{i| a_i > a_{i+1},1 \leq i \leq n-1\},\\[6pt]
&\mathrm{inv}(\omega)=\#\{(i,j)|a_i > a_j, 1 \leq i < j \leq n\},\\[6pt]
&\mathrm{maj}(\omega)=\sum_{a_i > a_{i+1},\atop{1 \leq i \leq n-1}} i,
 \end{align*}
 where $\#$ indicates the cardinality of a set. Writing  $\omega$  in the
 two-line form,  \begin{equation} \label{o2}
 \omega=\left(\begin{array} {cccc} x_1&x_2&\cdots&x_n\\
a_1&a_2&\cdots& a_n
 \end{array}\right),
 \end{equation}
 one  can define
 the excedance number $\mathrm{exc}(\omega)$ as follows
 $$
 \mathrm{exc}(\omega)=\#\{i| a_i > x_i,1 \leq i \leq n\}.
 $$
Usually, we say that $(a_i,a_{i+1})$ is a descent in $\omega$ if $a_i>a_{i+1}$ and $(a_i,x_i)$ is an excedance if $a_i>x_i$.
  %In Section \ref{psi},  we give an overview of Foata's first fundamental transformation  $\Phi_1$ and Han's bijection $\Psi$. We observe that these two bijections coincide when restricted to words on $\{1,2\}$. Moreover,
%  based on integer partition representation of a word and statistics on partitions,
%  we can use either bijection to describe Eulerian pairs $(F_n,\Phi_1^{-1}(F_n))$ and $(F_n^{'},\Phi_1^{-1}(F_n^{'}))$. In Section \ref{Gamma}, we extend  a bijection of
% Steingr\'{\i}msson on permutations to a map $\Gamma$ on words, which also maps the descent number to the excedance number. The resulting map is not a bijection on words. However, it turns out that $\Gamma$ is injective on $F^{'}_n$,
% and so we obtain another Eulerian pair $(F_n^{'},\Gamma(F_{n}^{'}))$.
%  Some remarks on  Euler-Mahonian pairs are given in Section \ref{remark}.

\section{Eulerian pairs derived from $\Phi_1^{-1}$}\label{psi}

In this section, we  construct Eulerian  pairs on Fibonacci words by using
Foata's first fundamental transformation \cite{Foata1965}.
 It is worth mentioning that Foata's first fundamental transformation $\Phi_1$  coincides with  Han's bijection  \cite{Han1994}
 when restricted to words on $\{1,2\}$. From now on, we shall still use $\Phi_1$ to denote Foata's first fundamental transformation (or Han's bijection) when restricted  to  $\{1,2\}^{*}$.

 Throughout this paper, by a binary
word we mean a word on $\{1,2\}$. Let %$\{1,2\}^{*}$ denote the set
%of binary words, and let
 $\{1,2\}_{n}^{*}$ denote the set of binary
words of length $n$.
 Clearly, a word
$\omega\in\{1,2\}^{*}$ with $d$ descents can be uniquely written as
\begin{align}\label{form}
\omega=1^{m_0}2^{n_0}1^{m_1}2^{n_1}\cdots 1^{m_d}2^{n_d},
\end{align}
where $m_0,n_d\geq 0$, and  $m_i,n_j>0$ for $1\leq i\leq d$ and $ 0\leq
j\leq d-1$.
It can be easily checked that $\Phi^{-1}_1(\omega)$ takes the following form
 \begin{equation}\label{form1}
\Phi_1^{-1}(\omega)=1^{m_0}21^{m_1-1}2\cdots
21^{m_d-1}2^{n_0-1}12^{n_1-1}\cdots 2^{n_{d-1}-1}12^{n_d}.
\end{equation}

The expression \eqref{form1} enables us to describe the Eulerian   pairs $(S,\Phi_1^{-1}(S))$  when $S=F_n$ and $S=F'_n$, where $F_n$ is the set of Fibonacci words of length $n$ and $F'_n$ is the set of Fibonacci words of length $n$ ending with $1$.
Analogous to the description of the Mahonian pairs obtained by
 Sagan and Savage \cite{Sagan-Savage2011},
we shall use the correspondence between binary words and integer partitions. Making use of this connection,  $\Phi_1^{-1}(F_n)$ and $\Phi_1^{-1}(F'_n)$ can be described in terms of statistics on integer partitions.

The following theorem gives
Eulerian pairs involving $F_n$ and $F'_n$, where we use
 $N_{\omega}(1)$  to denote the number of ones in a word $\omega$. For any partition $\lambda$, we denote by $l(\lambda)$ the number of parts of $\lambda$.
Recall that the Durfee square $D(\lambda)$ of $\lambda$ is  the square partition $(d^d)$, where $d$ is the largest integer $i\leq l(\lambda)$ such that
$\lambda_1\geq i, \ldots,\lambda_i\geq i$. Denote by $d(\lambda)$ the size $d$  of $D(\lambda)$,
and let  $B(\lambda)=(\lambda_{d+1},\ldots,\lambda_k)$.

\begin{thm}\label{reverseFnk} Let
$$ R_n = \{\omega \in \{1,2\}_n^*\ | \  \lambda=\lambda(\omega) , \lambda_1 \leq
n-N_{\omega}(1) , \, N_{\omega}(1)-1 \leq d(\lambda) \leq
N_{\omega}(1), B(\lambda)=\emptyset \},$$ and let
$$R_n'=\{\omega \in \{1,2\}_n^*\ | \  \lambda=\lambda(\omega),
\lambda_1= n-N_{\omega}(1) , N_{\omega}(1)-1 \leq d(\lambda) \leq
N_{\omega}(1),B(\lambda)=\emptyset \}.$$ Then
$(F_{n},R_n)$ and $(F'_{n},R'_n)$ are  Eulerian pairs.
\end{thm}

\noindent
 {\it
Proof.}  $\Phi_1$ is a bijection on words which maps the   excedance number  to the
 descent number, for more details, see \cite{Foata1965}. Thus for any set $S$,
$(S,\Phi_1^{-1}(S))$ is   an Eulerian pair.  So it suffices to show that
$R_n=\Phi_1^{-1}(F_{n})$ and $R'_n=\Phi_1^{-1}(F'_{n})$.

Suppose that $\omega=1^{m_0}2^{n_0}1^{m_1}2^{n_1}\cdots 1^{m_d}2^{n_d}\in F_n$,
where $m_0=0$ or $1$ and $m_{i}=1$ for $1\leq i\leq d$. Notice that $d=N_{\omega}(1)-m_0$. From (\ref{form1}) it follows that
\begin{equation}\label{phi}
\Phi_1^{-1}(\omega)=1^{m_0}2^{d+n_0-1}12^{n_1-1}\cdots
12^{n_{d-2}-1}12^{n_{d-1}-1}12^{n_d}.
\end{equation} Let
$\lambda=\lambda(\Phi_1^{-1}(\omega))$. From the correspondence between
binary words and partitions, we see that  $\lambda$ has exactly $d$
parts. Moreover, we have
 $$\lambda_1=n-N_{\omega}(1)-n_d \leq n-N_{\omega}(1)$$ and \[ \lambda_d=d+n_0-1 \geq d.\]
 Hence $B(\lambda)=\emptyset$ and $D(\lambda)=(d^d)$. It follows from \eqref{phi} that the size of the
 Durfee square of $\lambda$ is given by
\begin{align*}
d(\lambda)=\left\{
               \begin{array}{ll}
                 N_{\omega}(1)-1, & \hbox{if $m_0=1$;} \\[6pt]
                 N_{\omega}(1), & \hbox{if $m_0=0$.}
               \end{array}
             \right.
\end{align*}
So we see that $\Phi_1^{-1}(\omega)\in R_{n}$,
which yields that $\Phi_1^{-1}(F_n)\subseteq R_n$.
It is easy to see that the above  process is reversible and thus we arrive at the conclusion that $R_n= \Phi_1^{-1}(F_n)$.%have shown that

%Conversely,  let $\sigma=a_1a_2\cdots
%a_n \in R_{n}$. We wish to show that there is  a word $\rho\in F_n$ such that $\Phi_1^{-1}(\rho)=\sigma$. Let   $k=N_{\sigma}(1)$ and $\mu=\lambda(\sigma)$.
%By the definition of $R_n$, we have
% $\mu_1 \leq n-k $, $k-1 \leq d(\mu) \leq k$ and
%$B(\mu)=\emptyset$. So we may assume that  $\mu_1=n-k-t$ for some $t\geq 0$.
% By the construction of $\mu$, we see that
% $t$ is the largest integer $i$ such that $\sigma$ ends with $2^i$.  If $d(\mu)=k$,
%then we have $\mu_k \geq k$ and  $l(\mu)=k$. Hence $\sigma$ takes
% the form $2^ka_{k+1}\cdots a_{n-t-1}12^t$.
% Since $n-t-k=\mu_1\geq\mu_k \geq k$,  there exists a
%sequence of $k$ positive integers $n_0,n_1,\cdots,n_{k-1}$ such that
%$\sigma$ has the form $2^{k+n_0-1}12^{n_1-1}\cdots
%2^{n_{k-1}-1}12^t$. Let $\rho=2^{n_0}12^{n_1}1\cdots
%2^{n_{k-1}}12^t$. Obviously, $\rho \in F_{n}$. In view of
%(\ref{form1}), we find that $\Phi_1^{-1}(\rho)=\sigma$. For the case
%$d(\mu)=k-1$, by a similar argument
% it can be shown that there exists a  word $\rho'$ in $F_n$ such that $\Phi_1^{-1}(\rho')=\sigma$.
% So we have shown that $R_n\subseteq \Phi_1^{-1}(F_n)$. Consequently, we

 We now proceed to show that $R'_n=\Phi_1^{-1}(F'_{n})$.
Let $\omega$ be a binary word of length $n$.
In view of \eqref{form1}, we see  that $\omega$ ends with $1$ if and only if $\Phi_1^{-1}(\omega)$ ends with $1$. So we deduce that
\[
\Phi_{1}^{-1}(F'_n)=\{\omega\in \Phi_{1}^{-1}(F_n)\ | \  \omega{\rm\
ends\ with\ } 1\}.
\]
On the other hand, by the construction of the correspondence between binary words and partitions,
it can be checked  that  $\omega$  ends with $1$ if and only if $\lambda_1=n-N_{\omega}(1)$,
 where $\lambda=\lambda(\omega)$.  Since  $R_n= \Phi_1^{-1}(F_n)$,  we obtain that
\begin{align*}
\Phi_{1}^{-1}(F'_n)=&\{\omega\in R_n\ | \  \omega{\rm\ ends\ with\ } 1\}\\[6pt]
=&\{\omega\in R_n \ | \ \lambda=\lambda(\omega),\lambda_1=n-N_{\omega}(1)\},
\end{align*}
that is, $R'_{n}=\Phi_{1}^{-1}(F'_n)$.
This completes the proof.
\qed

\section{Eulerian pair derived from $\Gamma$ }\label{Gamma}

In this section, we  extend  the bijection of
Steingr\'{\i}msson  $\phi$ \cite{Steingrimsson1991} on permutations to
a map $\Gamma$ on words. While the extended map is not a bijection, it still
 transforms the descent number to the excedance number.
 As far as $F_n$ is concerned, the map $\Gamma$ is not injective, but it
  turns out to be injective on $F'_{n}$. Therefore, we obtain an Eulerian pair $(F'_{n},\Gamma(F'_{n}))$.

We begin with an overview of   Steingr\'{\i}msson's bijection $\phi$ on permutations.
Let $\pi=\pi_1\pi_2\cdots \pi_n$ be a permutation of $[n]$. %Note that any permutation of $[n]$  can be viewed as a one-to-one
% function on $[n]$.
 For notational convenience,   let $\phi(\pi)=f(1)f(2)\cdots f(n)$. %Instead of computing $f(k)$ directly, we proceed to find out $i_k$ such that $f(i_k)=\pi_k$ for $1\leq k\leq n$.
 Set $\pi_0=0$ and $\pi_{n+1}=n+1$.  \begin{itemize}
                                                                           \item [(1)] If there exists
an integer $m$ such that $k<m\leq n$ and
 $\pi_m<\pi_k $, then we set  $f(\pi_{k+1})=\pi_k$.
                                                                           \item [(2)]If $\pi_k>\pi_m$ for $k<m\leq n$, then we set
$f(\pi_{j+1})=\pi_k$, where $j$ is the largest number such that $\pi_j < \pi_k$.

                                                                         \end{itemize}

Steingr\'{\i}sson proved that the map $\phi$ is a bijection which maps the descent number to the excedance number.
\begin{prop}[\cite{Steingrimsson1991}, Remark 4.7]\label{s}Let $\pi$ be a permutation on $[n]$. Then for $1\leq k\leq n$,
$\pi_k>\pi_{k+1}$ if and only if $(\pi_k,\pi_{k+1})$ is an excedance  in $\phi(\pi)$.
 \end{prop}

  Steingr\'{\i}msson's bijection can be extended to a map $\Gamma$ on words. Recall that the standardization of a word  $\omega=a_1a_2\cdots a_n$
can be expressed as
$\pi=\beta_{\omega}(1)\beta_{\omega}(2)\cdots \beta_{\omega}(n)$ on
$[n]$, where $\beta_{\omega}(i)$ is given by
 \begin{equation}
 \beta_{\omega}(i)=\#\{j\ | 1\leq j\leq n, a_j<a_i\}+\#\{j\ |\ j\leq i,a_j=a_i\}.
 \end{equation}

Let $\omega=a_1a_2\cdots a_n$ be a word. The map $\Gamma$ is defined
as follows. Assume that
$\pi=\beta_{\omega}(1)\beta_{\omega}(2)\cdots \beta_{\omega}(n)$ is
the standardization of $\omega$. Let $\phi(\pi)=f(1)f(2)\cdots f(n)$.
   For $1\leq i\leq n$, there exists a unique integer $j_i$   such that $\beta_{\omega}(j_i)=f(i)$.
  Then  $\Gamma(\omega)$ is defined to be
  the word $a_{j_1}a_{j_2}\cdots a_{j_n}$. For example, let
 $\omega=132232131$. Then the standardization of $\omega$ is $\pi=174586293$ and
  $\phi(\pi)=169748253$.   So we have $\Gamma(\omega)=123323121.$

The following theorem  shows that the map $\Gamma$ also transforms the descent number to the excedance number.

\begin{thm}
For any word $\omega$, we have
$$\mathrm{des}(\omega)=\mathrm{exc}(\Gamma(\omega)).$$
\end{thm}

\noindent{\it Proof.}
Assume  that  $\omega=a_1a_2\cdots a_n$ is a word. Let $\pi=\sigma_1\sigma_2\cdots\sigma_n$ be the standardization of $\omega$.
It is obvious that $(a_i,a_{i+1})$ is a descent in $\omega$ if and only if $(\sigma_i,\sigma_{i+1})$ is a descent in $\pi$. By Proposition \ref{s}, we see that $(\sigma_i,\sigma_{i+1})$ is a descent in $\pi$ if and only if $(\sigma_i,\sigma_{i+1})$ forms an excedance in $\phi(\pi)$. With the aid of the construction of $\Gamma$, it can be seen that $(\sigma_i,\sigma_{i+1})$ forms an excedance in $\phi(\pi)$ if and only if $(a_i,a_{i+1})$ is an excedance in $\Gamma(\omega)$. Thus, we have  $\mathrm{des}(\omega)=\mathrm{exc}(\Gamma(\omega))$. This completes the proof.\qed

 Next we consider the restriction of $\Gamma$ to words on  $\{1,2\}$. In this case, it is easy to verify
 that $\Gamma(\omega2^m)=\Gamma(\omega)2^m$ for  $m \geq 1$. The following lemma shows how to compute $\Gamma(\omega 1^m)$ based on $\Gamma(\omega)$.

\begin{lem}\label{lemma}
 Suppose that $\omega$ is a binary word of length $n$  that contains
 $k$ ones. Let  $\Gamma(\omega)=b_1b_2\cdots b_n$. Assume that
  $t$ is the largest integer $i$ such that $\omega$ ends with $2^{i}$.
   Set $U=b_1b_2\cdots b_k$ and $V=b_{k+1}b_{k+2}\cdots b_{n-t}$. Then we have the following
recurrence relations:
\begin{itemize}
\item[(1)] If $t=0$, then $\Gamma(\omega1)=U1V$. In general, if  $t=0$, then
$\Gamma(\omega1^m)=U1^{m}V$ for any $m \geq 1$;
\item[(2)] If $t>0$, then $\Gamma(\omega1)=U2V12^{t-1}$. In general, if $t>0$,  then we
have  $\Gamma(\omega1^m)=U21^{m-1}V12^{t-1}$ for any $m \geq 1$.
\end{itemize}
\end{lem}

\noindent{\it Proof.}  Let $\omega=a_1a_2\cdots a_n$ and   $a_{n+1}=1$. Suppose that $\Gamma(\omega1)=c_1c_2\cdots c_{n+1}$.  To determine $\Gamma(\omega1)$,  we consider
 occurrences of ones in $\Gamma(\omega1)$.  Assume that
  $a_{s_1}, a_{s_2}, \ldots, a_{s_k}$ are the ones in $\omega$,
  where $s_1< s_2 < \cdots < s_k$.     Let us  define
$\beta(i)=\beta_{\omega}(i)$ and $\beta'(j)=\beta_{\omega1}({j})$
for $1\leq i\leq n$ and $1\leq j\leq n+1$. It can be seen that
$\beta'(n+1)=k+1$ and for $i\leq n$,
\[
\beta'(i)=\left\{
        \begin{array}{ll}
          \beta(i), & \hbox{if $a_i=1$;} \\[6pt]
          \beta(i)+1, & \hbox{otherwise.}
        \end{array}
      \right.
\]
Thus we have
\[ \{\beta'(s_1)<\beta'(s_2)<\cdots<\beta'(s_k)\}=\{1,2,\ldots,k\}\]
and \[ \{\beta'(i)| 1\leq i\leq n, a_i=2\}=\{k+2,\cdots,n+1\}.\]
By the construction of $\Gamma$, it is not hard to see that
$b_{\beta(s_i+1)}=a_{s_{(i+1)}}=1$ and
$c_{\beta'(s_i+1)}=a_{s_{(i+1)}}=1$ for $0 \leq i \leq k-1$, where
$s_0=0$. For $0\leq i\leq k-1$,  it is clear that $\beta'(s_i+1)\leq k$ if and only if $\beta'(s_i+1)=\beta(s_i+1)$.
This means that the ones in $c_{1}c_{2}\cdots c_{k}$ appear in
the same positions as in $U$.
 Moreover,   for the case $\beta'(s_i+1)\geq k+2$, we see
that $\beta'(s_i+1)=\beta(s_i+1)+1$.   In other words, a one appearing in the $j$-th position  in
$V$ corresponds to a one in the $j$-th position  in $c_{k+2}c_{k+3}\cdots c_{n-t+1}$.

Let us further consider the position of $a_{n+1}$
in $\Gamma(\omega1)$. Observe that $s_k=n-t$. By the construction of
$\Gamma$, we find that $c_{\beta'(n-t+1)}=a_{n+1}=1$. If $t=0$, then
$c_{k+1}=c_{\beta'(n+1)}=a_{n+1}$, which means that  $a_{n+1}$
is in the $(k+1)$-th position in $\Gamma(\omega1)$. When $t
>0$,  since $\beta'(n-t+1)=n-t+2$,  we find that $c_{n-t+2}=c_{\beta'(n-t+1)}=a_{n+1}$. Thus $a_{n+1}$ is in the
$(n-t+2)$-th position in $\Gamma(\omega1)$. In summary, we deduce that
\begin{align}\label{the case 1}
\Gamma(\omega1)=\left\{
\begin{array}{ll}
U1V,& {\rm if }\ t=0;\\[6pt]
U2V12^{t-1},& {\rm if }\ t >0.
\end{array}
\right.
\end{align}
So the lemma holds for $m=1$. By iterating the above process,
it can be seen that the lemma holds for $m>1$.  This completes the proof. \qed

By Lemma \ref{lemma}, for any word $\omega$ in form \eqref{form},  $\Gamma(\omega)$  is of the following form
\begin{equation}\label{formula of gamma}
\Gamma(\omega)=1^{m_0}21^{m_1-1}\cdots 21^{m_{d-1}-1}
21^{m_d}2^{n_0-1}12^{n_1-1}\cdots 2^{n_{d-2}-1}12^{n_{d-1}-1+n_d}.
\end{equation}

The following theorem gives a description of $\Gamma(F'_n)$.

\begin{thm}\label{gamma1}Let
 $$T_n=\{\omega \in
\{1,2\}_n^*\ | \ \lambda=\lambda(\omega), \lambda_1 \leq
n-N_{\omega}(1) , \, N_{\omega}(1)-1 \leq
l(\lambda)=\lambda_{l(\lambda)} \leq N_{\omega}(1) \}.$$ Then we
have $\Gamma(F'_n)=T_n$. Moreover, $(F'_{n},T_n)$ is an
Eulerian pair.
\end{thm}
\noindent{\it Proof.} Using the argument in the proof of Theorem \ref{reverseFnk},
it can be shown that $\Gamma(F'_n)=T_n$.
To prove  $(F'_{n},T_n)$ is an
Eulerian pair, it suffices to verify that $\Gamma$ is injective on $F'_n$.
Assume  that $\omega=1^{m_0}2^{n_0}12^{n_1}\cdots
12^{n_{d-2}}12^{n_{d-1}}1$ and $\omega'=1^{m'_0}2^{n'_0}12^{n'_1}\cdots
12^{n'_{d'-2}}12^{n'_{d'-1}}1$
 are two words in $F'_n$ such that
$\Gamma(\omega)=\Gamma(\omega')$. It follows from \eqref{formula of gamma} that $\Gamma(\omega)=1^{m_0}2^d12^{n_0-1}\cdots 2^{n_{d-2}-1}1$ and
$\Gamma(\omega')=1^{m'_0}2^{d'}12^{n'_0-1}\cdots 2^{n'_{d'-2}-1}1$.
So we have $d=d'$, $m_0=m'_0$ and $n_i=n'_i$ for any $0 \leq i
\leq d-1$. This implies that $\omega=\omega'$. Hence
  $\Gamma$ is injective on
$F'_{n}$. This completes the proof.\qed

It should be noted that $\Gamma$ is neither  surjective nor injective on $F_n$. For example,
  there is no  $\omega$ satisfying
   $\Gamma(\omega)=2121$.  On the other hand, we have
$$\Gamma(2^212^212^31)=\Gamma(2^212^212^212)=\Gamma(2^212^21212^2)=2^3121212^2.$$

We conclude  this section with a remark  that $\Gamma(F_{n})=\Gamma(F'_{n})$. In fact,  for any word $\omega=1^{m_0}2^{n_0}12^{n_1}\cdots 12^{n_{d-1}}12^{n_{d}}\in F_{n}$, let $\sigma=1^{m_0}2^{n_0}12^{n_1}\cdots
12^{n_{d-1}+n_{d}}1$ in $F'_{n}$. Then we have
$\Gamma(\omega)=\Gamma(\sigma)$.

\section{Concluding Remarks}\label{remark}

In this section, we make some remarks on Euler-Mahonian pairs on binary words, which
are related to the bijections $\Phi_1$, $\Phi_2$ and $\Gamma$.

For any word $\omega=1^{m_0}2^{n_0}\cdots 1^{m_d}2^{n_d}$, Sagan and Savage have shown  that
\begin{align}\label{Phi_2}
\Phi_2(\omega)=1^{m_d-1}21^{m_{d-1}-1}2\cdots 1^{m_1-1}21^{m_0}2^{n_0-1}12^{n_1-1}1\cdots 2^{n_{d-1}-1}12^{n_d}.
\end{align}
%where $\Phi_2$ is the second fundamental transformation.
It is clear from \eqref{Phi_2} that
$
{\rm des}(\omega)={\rm exc}(\Phi_2(\omega))
$. So we deduce that  the  Mahonian pairs $(S,T)$ given
 by Sagan and Savage \cite{Sagan-Savage2011} are  Euler-Mahonian pairs in the sense that
\[
\sum_{\omega \in
S}p^{\mathrm{des}(\omega)}q^{\mathrm{maj}(\omega)}=\sum_{\omega \in
T}p^{\mathrm{exc}(\omega)}q^{\mathrm{inv}(\omega)}.
\]

% Eulerian pairs considered in the previous sections are %different from those derived from $\Phi_2$,
 It should be noted  that in general $\Phi_2(F_n)\neq \Phi_1^{-1}(F_n)$, $\Phi_2(F'_n)\neq \Phi_1^{-1}(F'_n)$ and $\Phi_2(F'_n)\neq \Gamma(F'_n)$.
However, there exists a set $G_n$ such that $(G_n,\Phi_{1}^{-1}(G_n))$, $(G_n,\Phi_2(G_n))$ and $(G_n,\Gamma(G_n))$ are the same Eulerian pairs. Meanwhile, we find a set $H$ of binary words for which $\Phi_{1}^{-1}=\Phi_2$.

\begin{thm}
Let $G_{n}$ be  the set of  words
in $\{1,2\}^{*}_{n}$ with no consecutive twos and let
\[
H=\{\omega=1^{m_0}2^{n_0}1^{m_1}2^{n_1}\cdots 1^{m_d}2^{n_d}|\
m_0=m_d-1, m_i=m_{d-i}\ {\rm for}\  1 \leq i \leq
    d-1\}.
\]
Then we have
$\Phi_2(G_n)=\Phi_1^{-1}(G_n)=\Gamma(G_n)$ and $\Phi_{1}^{-1}(\omega)=\Phi_2(\omega)$ for any $\omega\in H$.
\end{thm}

\noindent{\it Proof.}  Given a word $\omega\in G_n$ with $d$ descents, it can be written uniquely as
\[
1^{m_0}21^{m_1}2\cdots 1^{m_d}2^{n_d},
\]
where $m_0\geq 0$, $n_d=0$  or $1$, and $m_i>0$ for $1\leq i\leq d$.  By \eqref{form1} and \eqref{formula of gamma}, we find   that $\Phi_1^{-1}(\omega)=\Gamma(\omega)$ for all $\omega\in G_n$. Therefore, we have $\Phi_1^{-1}(G_n)=\Gamma(G_n)$.    To show    that $\Phi_1^{-1}(G_n)=\Phi_2(G_n)$, we define a map $\varphi$ on binary words
\begin{align*}
\varphi(1^{m_0}2^{n_0}1^{m_1}2^{n_1}\cdots 1^{m_{d-1}}2^{n_{d-1}}
1^{m_d}2^{n_d})=1^{m_{d}-1}2^{n_0}1^{m_{d-1}}2^{n_1}\cdots
2^{n_{d-2}}1^{m_1}2^{n_{d-1}}1^{m_0+1}2^{n_d}.
\end{align*}
It is easy to check that $\varphi$ is an involution on $\{1,2\}^{*}$.
Observing that $\varphi(G_n)=G_n$,
  by \eqref{Phi_2} and \eqref{form1}, we obtain that
$\Phi_2(\omega)=\Phi_1^{-1}(\varphi(\omega))$
for any $\omega\in G_n$. Thus we have $\Phi_1^{-1}(G_n)=\Phi_2(G_n)$.

By the definition of $\varphi$, we find
 that $H=\{\omega\in\{1,2\}^{*}\ |\ \varphi(\omega)=\omega\}$.
 Since  $\Phi_2(\omega)=\Phi_1^{-1}(\varphi(\omega))$ for any binary word $\omega$, we conclude that $\Phi_{1}^{-1}(\omega)=\Phi_2(\omega)$ for any $\omega\in H$.
This completes the proof. \qed

 \noindent{\bf Acknowledgments.} We wish to thank the referee
 for helpful suggestions. This work was supported by the 973
Project, the PCSIRT Project of the Ministry of Education, and the
National Science Foundation of China.

\end{document}